%% file: Arxiv_SOLID/main.tex
\documentclass[11pt,a4paper]{article}
\pdfoutput=1

\usepackage[numbers,sort&compress]{natbib}
\usepackage{sjtu_iic_arxiv}
\usepackage{mathtools}
\usepackage{booktabs}
\usepackage{multirow}
\usepackage{algorithm}
\usepackage{algorithmic}
\usepackage{array}
\usepackage{colortbl}
\usepackage{placeins}
\usepackage{subcaption}
\usepackage{multicol}

\graphicspath{{figs/}}
\newcommand{\tblfont}{\small}
\newcommand{\method}{\textsc{SOLID}}
\newcommand{\Rvote}{\hat z_{\mathrm{vote}}}

\hypersetup{
  pdftitle={Beyond Verified Answers: Solver-Informed Self-Distillation for Bootstrapping Operations Research Language Models},
  pdfauthor={First Author, Second Author}
}

\setheadertext{%
  \raisebox{-0.55cm}{\includegraphics[height=1.3cm]{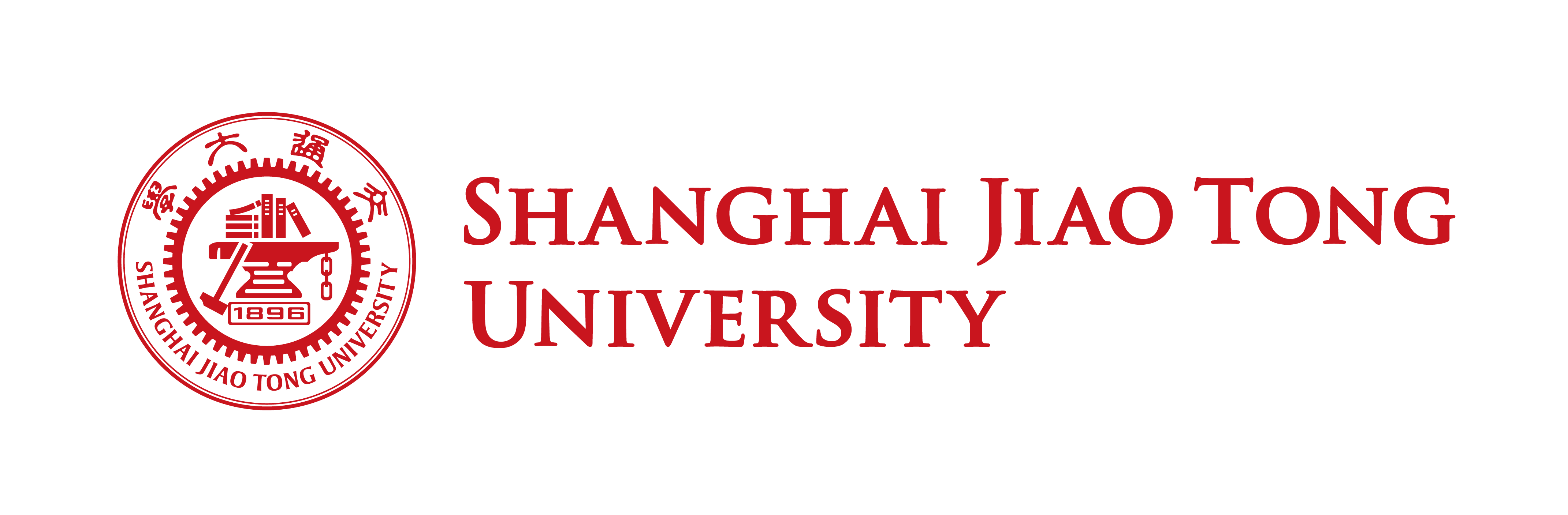}}%
  \hspace{0.15cm}%
  \raisebox{-0.29cm}{\includegraphics[height=0.72cm]{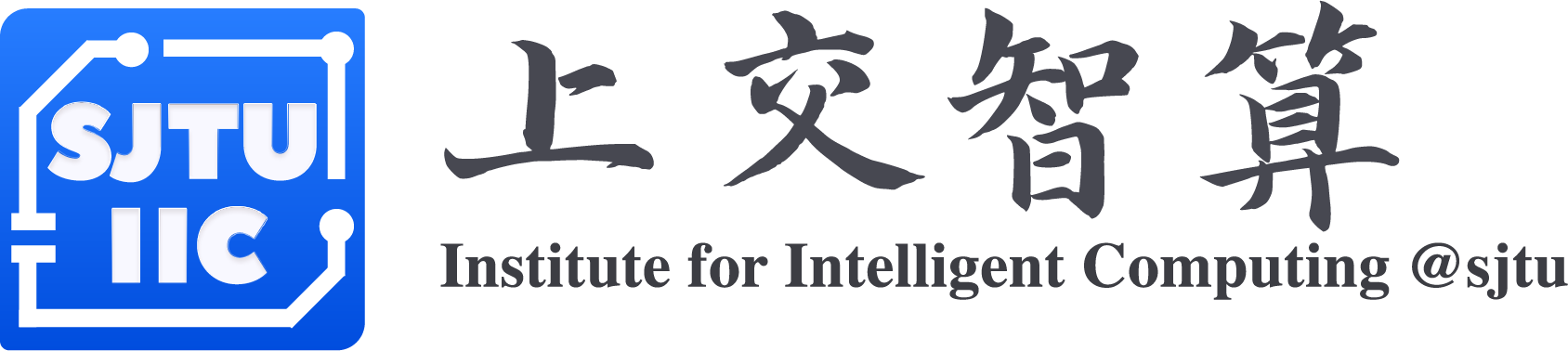}}%
}
\setheadertitle{Solver-Informed Self-Distillation for OR Language Models}
\githublink{https://github.com/AIOR-Research/SOLID}
\huggingfacelink{https://huggingface.co/AIOR-Research/SOLID-Qwen3}
\correspondingemail{%
  \texttt{\{jameszhu, chenyuzhou, linjianghao, ddge\}@sjtu.edu.cn}\quad
  \texttt{mlcao25@fudan.edu.cn}\\[0.25cm]
  \textsuperscript{*}Equal contribution\quad
  \textsuperscript{\dag}Corresponding author
}

\title{%
  Beyond Verified Answers: Solver-Informed Self-Distillation for
  Bootstrapping Operations Research Language Models
}

\author{%
  Rui Zhu\textsuperscript{1,*}\hspace{0.55em}
  Minglong Cao\textsuperscript{2,*}\hspace{0.55em}
  Chenyu Zhou\textsuperscript{1,*}\hspace{0.55em}
  Jianghao Lin\textsuperscript{1,\dag}\hspace{0.55em}
  Dongdong Ge\textsuperscript{1}\hspace{0.55em}\\[2pt]
  {\sjtuiicAffilFont
  \textsuperscript{1}Shanghai Jiao Tong University, Shanghai, China\\
  \textsuperscript{2}Fudan University, Shanghai, China}
}

\begin{document}

\begin{abstract}
\input{sections/abstract}
\end{abstract}

\maketitle

\input{sections/introduction}
\input{sections/related_work}

\input{sections/preliminaries}
\input{sections/method}
\input{sections/theory}
\input{sections/experiments}
\input{sections/conclusion}

\small
\bibliographystyle{plainnat}
\bibliography{refs}

\input{sections/appendix}

\end{document}

%% file: sections/abstract.tex
Modern large language models (LLMs) can translate natural-language descriptions into operations research (OR) formulations. Post-training techniques including reinforcement learning and on-policy self-distillation have further improved this capability. However, three limitations remain in training LLMs for OR formulations. 
First, training commonly relies on synthetic formulations validated by human experts or stronger models, constraining scalable supervision. 
Second, credit assignment is either coarse or costly: outcome rewards score an entire trajectory without locating the responsible modeling decision, whereas process-level supervision requires an additional evaluator. 
Third, privileged self-distillation can induce style mismatch by using solver context unavailable at deployment.
We find that a model can improve from solver-artifact feedback generated by its own rollouts, making self-distillation a practical, evaluator-free source of dense supervision. Therefore, we propose \textbf{\method{}}---\textbf{S}olver-Informed \textbf{O}n-Policy \textbf{L}earn\textbf{I}ng through Self-\textbf{D}istillation, a novel framework for self-improving OR language models without verified answers or external evaluators. \method{} executes candidate programs from multiple rollouts, clusters their objectives, and selects a majority-group artifact as a pseudo-reference. The model then performs updates using group-relative advantages and dense self-supervision signals. Across multiple OR benchmarks, \method{} improves solution accuracy for both general-purpose and OR-tuned models over outcome-only group-relative training. These results show that solver artifacts can support scalable self-improvement without trusted answers. 
%Our code is available at \href{https://github.com/jameszhu-pixel-time/SOLID}{\texttt{https://github.com/jameszhu-pixel-time/SOLID}}.

%% file: sections/introduction.tex
%\begin{figure*}[t]
%\centering
%\includegraphics[width=0.75\textwidth]{solid_method.png}
%\caption{The \method{} pipeline. Solver execution converts on-policy programs into a voted pseudo-objective and a selected solver artifact. The same policy uses this privileged context as a teacher, while the deployable student sees only the original problem. Group-relative rewards retain exploration, and a stage-aware mask restricts dense self-distillation to task-relevant modeling and code stages.}
%\label{fig:solid-pipeline}
%\end{figure*}

%\begin{figure}[t]
%    \centering
%    \includegraphics[width=\columnwidth,]{figs/0728_Intro_v1.png}
%    \caption{
%    Motivation of \method{}.
%Scalable supervision for OR language models is challenging because obtaining verified formulations and step-level annotations requires substantial domain expertise.
%Existing solver-based training signals further suffer from two limitations:
%outcome-only rewards provide coarse sequence-level feedback and cannot identify errors in intermediate modeling components,
%while privileged self-distillation introduces information leakage by relying on solver information unavailable to the deployed model.
%    }
%    \label{fig:solid_overview}
%\end{figure}

\section{Introduction}

Large language models (LLMs) can derive mathematical formulations and
executable solver programs from natural-language operations research (OR) problem descriptions. Inference-time search and reinforcement learning (RL) have improved this capability \citep{wang2025bppsearch,ramamonjison2022nl4opt,huang2025orlm,lu2025optmath,zhou2026steporlm}.
Reliable OR modeling remains challenging, however, because the variables, objective, constraints, and code must jointly represent the intended decision problem \citep{jiang2025llmopt}. 
Despite this progress, these RL related paradigms still face three key bottlenecks in the OR setting.

The first obstacle is \textbf{the cost of scalable supervision}. Existing
pipelines expand curated seed formulations into synthetic training examples
and validate them with trusted solutions, stronger models, or learned evaluators
\citep{huang2025orlm,lu2025optmath,zhou2026steporlm}. In addition, validation still requires
OR expertise because equivalent formulations may differ substantially, while
an incorrect formulation may execute and return a plausible objective. This
dependence limits learning from problems for which descriptions and solver
access are available but verified answers are not.

The second obstacle is \textbf{accurate credit assignment within a single trajectory}. This is challenging in OR, since reasoning typically spans variable definition, objective formulation, constraint construction, and code generation.
\begin{figure}[H]
    \centering
    \includegraphics[width=0.85\textwidth]{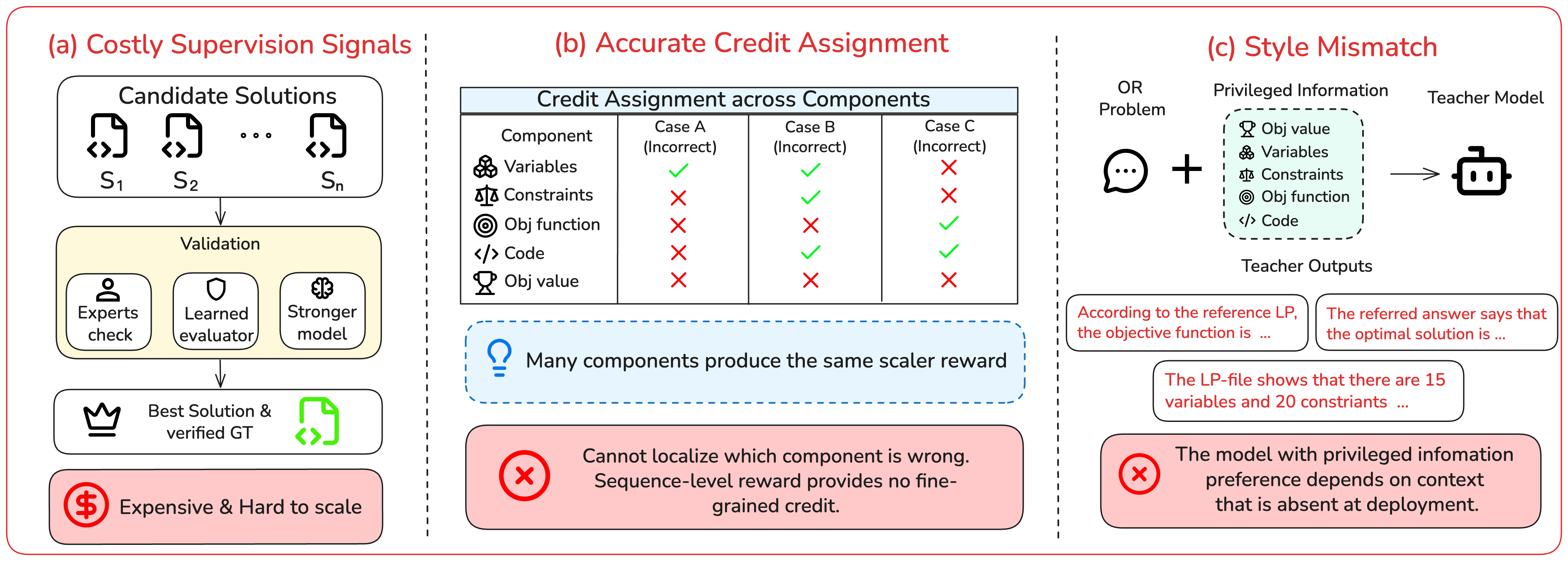}
    \captionsetup{skip=4pt}
    \caption{Three challenges motivating \method{}: costly supervision,
    coarse sequence-level credit assignment, and privileged-context style
    mismatch.}
    \label{fig:solid_overview}
\end{figure}
OR research has refined sequence-level signal by
evaluating formulation components separately, deriving verifiable feedback from generated code and solver artifacts, and learning process reward models for intermediate reasoning
\citep{ramamonjison2022nl4opt,NEURIPS2025_98555b92,zhou2026steporlm}. These approaches provide finer credit, but require component annotations, predefined verification signals, or another learned evaluator. On-policy self-distillation (OPSD) provides an economical alternative for assigning dense credit: the same model re-evaluates its trajectory under privileged information, without a stronger teacher
\citep{zhao2026selfdistilled,hubotter2026sdpo}.

The third obstacle is \textbf{style mismatch}, a form of information leakage in which the contextual policy conditions its predictions on evidence that is unavailable at deployment time \citep{pan2026rlcsd,kim2026doesselfdistillationsometimesdegrade}. This issue is particularly pronounced in operations research. When a correct solver artifact is provided as privileged information, the privileged policy may rename variables, omit intermediate derivations, or simply restate the artifact itself. Consequently, directly matching the privileged policy’s distribution conflates superficial changes in form and degrading performance.

To address all obstacles, we propose \method{} (Solver-Informed On-Policy
LearnIng through Self-Distillation), which converts solver information from
the same rollout group into both label-free sequence supervision and  token-level supervision. Solver-outcome supplies the
group-relative learning signal, while one majority-group solver artifact
provides privileged context for self-distillation. A sample
mask restricts self-distillation only to structurally mismatched
variable, objective, constraint, and code sections. This closed loop requires neither verified
answers nor an external evaluator and mitigates style
mismatch.

Our contributions are threefold:
\begin{itemize}
    \item We formulate OR language-model post-training \emph{without} verified
    answers as a closed self-bootstrapping loop over executable on-policy
    rollouts and solver-outcome consensus, requiring no external evaluator.
    \item We introduce complementary sequence-level and token-level supervision from the same rollout group. Solver-outcome agreement drives
group-relative learning, while self-distillation with an
LP-structured mask blocks potentially harmful updates for
plausible or unrelated reasoning steps.
    \item \method{} consistently outperforms self-distillation and label-free reinforcement learning methods across challenging benchmarks. Its effectiveness on base and post-trained models further demonstrates its broad applicability and potential for scalable OR post-training without verified answers.
\end{itemize}

%\begin{figure}[H]
%    \centering
%    \includegraphics[width=\textwidth]{figs/0729_Intro_horizontal.png}
    
%    \caption{
%Motivation for \method{}.
%(a) Verified supervision requires costly expert- or model-based validation.
%(b) Response-level solver outcomes cannot identify which modeling component is responsible for success or failure, resulting in ambiguous credit assignment.
%(c) Solver-privileged self-distillation may introduce style mismatch through preferences induced by context unavailable at deployment. These limitations motivate scalable, fine-grained, and deployment-consistent supervision.
%}
%    \label{fig:solid_overview}
%\end{figure}

%\begin{figure}[!t]
%    \centering
    %\includegraphics[width=\columnwidth,]{figs/0728_Intro_merge.png}
%    \includegraphics[
%        width=\columnwidth,
        %height=0.80\textheight,
        %keepaspectratio
%        ]{figs/0729_Intro_v1.png}
%    \caption{
%    Motivation for \method{}.
%(a) Obtaining verified formulations requires costly validation and is difficult to scale.
%(b) On-policy self-distillation introduces \emph{style mismatch}, because the teacher conditions on context unavailable to the deployed policy and may prefer reference-dependent surface forms.
%These limitations motivate scalable and deployment-consistent supervision without verified ground truth.
%}
%    \label{fig:solid_overview}
%\end{figure}

%% file: sections/related_work.tex
\section{Related Work}

\paragraph{Language Models for Optimization Modeling.}

Early auto-formulation systems cast natural-language descriptions into structured linear-programming (LP) components and emphasized constrained generation and human validation~\citep{ramamonjison-etal-2022-augmenting,ramamonjison2022nl4opt}. Subsequent work broadened the target to end-to-end formulation, code generation, and solving. Agentic approaches organize these activities through external tools and interaction protocols, within a broader landscape of memory, skills, and harness engineering~\citep{zhou2026externalization}.OptiMUS uses modular decomposition and iterative solver-based correction, while BPP-Search explores alternative formulations through tree search~\citep{ahmaditeshnizi2026optimus03usinglargelanguage,wang2025bppsearch}.
Agora-Opt combines solver-outcome-based decentralized debate with persistent memory to support training-free improvement in optimization modeling~\citep{lin2026soliloquy}. OR-Space extends evaluation to model construction, revision, and grounded explanation in persistent industrial workspaces~\citep{zhou2026orspace}.

%Early auto-formulation systems cast natural-language descriptions into structured linear-programming (LP) components and emphasized constrained generation and human validation \citep{ramamonjison-etal-2022-augmenting,ramamonjison2022nl4opt}. Subsequent work broadened the target to end-to-end formulation, code generation, and solving. OptiMUS uses modular decomposition and iterative solver-based correction, while BPP-Search explores alternative formulations through tree search \citep{ahmaditeshnizi2026optimus03usinglargelanguage,wang2025bppsearch}. These systems establish the value of executable intermediate representations, but primarily improve inference or train on curated formulation data rather than learn from unlabeled rollout groups.

\paragraph{Solver-Informed Supervision.}
Solver-informed supervision uses executable programs and optimization
artifacts to filter training data or construct learning signals. Existing OR
post-training pipelines typically build instruction data from synthetic formulations and develop trusted answers from solver information
\citep{huang2024mamo,huang2025orlm,lu2025optmath}. More recent methods move solver information into the learning signal: SIRL derives verifiable rewards from
generated code and LP artifacts; StepORLM adds a co-evolved generative process
reward model; and MURKA combines semantic and execution rewards with an
external teacher and a checker agent
\citep{NEURIPS2025_98555b92,zhou2026steporlm,xie2025murka}. Although these
approaches make supervision more task-grounded, their supervision signals still rely on trusted sources or external advanced models. Thus, how to derive dense supervision from solver feedback without verified answers remains an open question.

\paragraph{Self-Supervised Post-Training.}
Self-supervised post-training can derive learning signals from a model's own
on-policy rollouts. At sequence level, Test-time Reinforcement Learning (TTRL) turns
majority agreement into a group-relative reward for training without ground truth~\citep{zuo2025ttrl,shao2024deepseekmath}. At the token level, on-policy self-distillation (OPSD) re-scores a rollout using
the same policy conditioned on privileged information,
providing dense guidance without a stronger teacher~\citep{hubotter2026sdpo,he2026sdzero,zhao2026selfdistilled}. However, this contextual signal is not always reliable: privileged information may alter the notation, ordering, or reasoning style, causing self-distillation to introduce style mismatch and degraded performance~\citep{pan2026rlcsd,kaur2026rethinking}. \method{} connects these two branches by retaining TTRL’s group-relative outcome learning while applying self-distillation only to the segments of the reasoning traces where errors occur, thereby mitigating the \textbf{style mismatch} problem.

%% file: sections/preliminaries.tex
\section{Preliminaries}
\label{sec:preliminary}
Let $x$ be an unlabeled OR problem and
$y=(y_1,\ldots,y_T)$ a structured response ending in executable solver
code. Following recent LLM-based OR reasoning pipelines
\citep{NEURIPS2025_98555b92,zhou2026steporlm}, we view $y$ as a progression
from problem interpretation to decision variables, objective, constraints,
mathematical model, and solver code. The rollout policy samples
$\{y_i\}_{i=1}^{N}\sim\pi_{\mathrm{old}}(\cdot\mid x)$. An external Python
executor then runs each generated program and returns an invalid status or an optimal objective $\hat z_i$ together with its LP artifact. A simplified LP
artifact has the following solver-readable form:
\begin{center}
\fbox{\begin{minipage}{0.88\columnwidth}
\small\ttfamily
Minimize\\
\hspace*{1em}50 x + 75 y + 20 z\\
Subject To\\
\hspace*{1em}capacity: 50 x + 75 y + 20 z <= 20000\\
\hspace*{1em}balance: -2.33 x - y + z >= 0\\
Bounds\\
\hspace*{1em}x >= 50\\
Generals\\
\hspace*{1em}x y z\\
End
\end{minipage}}
\end{center}

\paragraph{Step-level credit assignment.}
Step-level methods refine trajectory-level supervision by assigning feedback
to intermediate decisions. StepOPSD is one example. It decomposes agent
trajectories into action-centered segments and redistributes
hindsight-conditioned distillation signals at this granularity
\citep{zhang2026stepopsd}. OR modeling provides a natural step structure
because a formulation is assembled through decisions
about variables, objectives, constraints, mathematical models, and executable
code
\citep{ramamonjison2022nl4opt,huang2025orlm}. StepORLM makes this structure
explicit by applying generative process supervision to intermediate OR
modeling steps \citep{zhou2026steporlm}. This structure motivates localizing
credit at the step level.

\paragraph{Reinforcement learning without labels.}Following common reinforcement-learning practices in mathematical reasoning
and OR, we adopt Test-time Reinforcement Learning (TTRL), an algorithm for reinforcement learning without labels, as our
baseline for unlabeled OR problems. Its implementation is based on Group Relative Policy Optimization (GRPO).
GRPO normalizes group rewards and applies a PPO-style clipped update without a
learned value model \citep{schulman2017ppo,shao2024deepseekmath}. Following
TTRL \citep{zuo2025ttrl}, solver-outcome agreement supplies the vote reward,
while format and execution rewards enforce structured output and solver
validity \citep{NEURIPS2025_98555b92}:
\[
 r_i=r_i^{\mathrm{fmt}}+r_i^{\mathrm{exec}}+r_i^{\mathrm{vote}},
 \qquad
 A_i=\frac{r_i-\mu_r}{\sigma_r+10^{-8}},
\]
and the corresponding clipped loss without Kullback–Leibler (KL) divergence is
\[
\begin{aligned}
\mathcal L_{\mathrm{GRPO}}
&=-\mathbb E_{i,t}\!
  \left[\min\{\rho_{i,t}A_i,\bar\rho_{i,t}A_i\}\right],\\
\bar\rho_{i,t}&=\operatorname{clip}(\rho_{i,t},1-\epsilon,1+\epsilon),\\
\rho_{i,t}&=
\frac{\pi_\theta(y_{i,t}\mid x,y_{i,<t})}
{\pi_{\mathrm{old}}(y_{i,t}\mid x,y_{i,<t})}.
\end{aligned}
\]
Every token shares $A_i$, so this \textbf{sequence-level supervision} cannot
locate the responsible modeling choice.

\paragraph{On-policy self-distillation.}
For the same rollout $y_i$, let $\widetilde x_i$ denote a privileged-context
version of $x$. OPSD uses the rollout policy as a fixed contextual teacher
while the updated policy remains the deployable student~\citep{zhao2026selfdistilled}:
\[
\begin{aligned}
\ell^T_{i,t}
&=\log\pi_{\mathrm{old}}
(y_{i,t}\mid\widetilde x_i,y_{i,<t}),\\
\ell^\theta_{i,t}
&=\log\pi_\theta(y_{i,t}\mid x,y_{i,<t}),\qquad
\Delta_{i,t}=\ell^T_{i,t}-\ell^\theta_{i,t}.
\end{aligned}
\]
Using the same sampled reverse-KL estimator as \method{}
\citep{schulman2020approximatingkl}, whole-response OPSD minimizes
\[
\mathcal L_{\mathrm{OPSD}}
=
\mathbb E_{i,t}\!\left[
\phi_{k_3}(\Delta_{i,t})
\right],
\qquad
\phi_{k_3}(\Delta)=e^\Delta-\Delta-1.
\]
Thus, OPSD provides dense token-level supervision but applies it uniformly to
the response and does not include the group-relative objective. In OR, the
majority-LP pseudo-reference can supply $\widetilde x_i$ and expose the model's
correction preferences without another evaluator. Because uniform re-scoring
suffers from \textbf{style mismatch}, \method{} uses the binary LP-structured
mask to localize supervision signals~\citep{pan2026rlcsd}.

%% file: sections/method.tex
\section{Methodology}

\label{sec:methodology}

\begin{figure}[H]
    \centering
    \includegraphics[width=\textwidth]{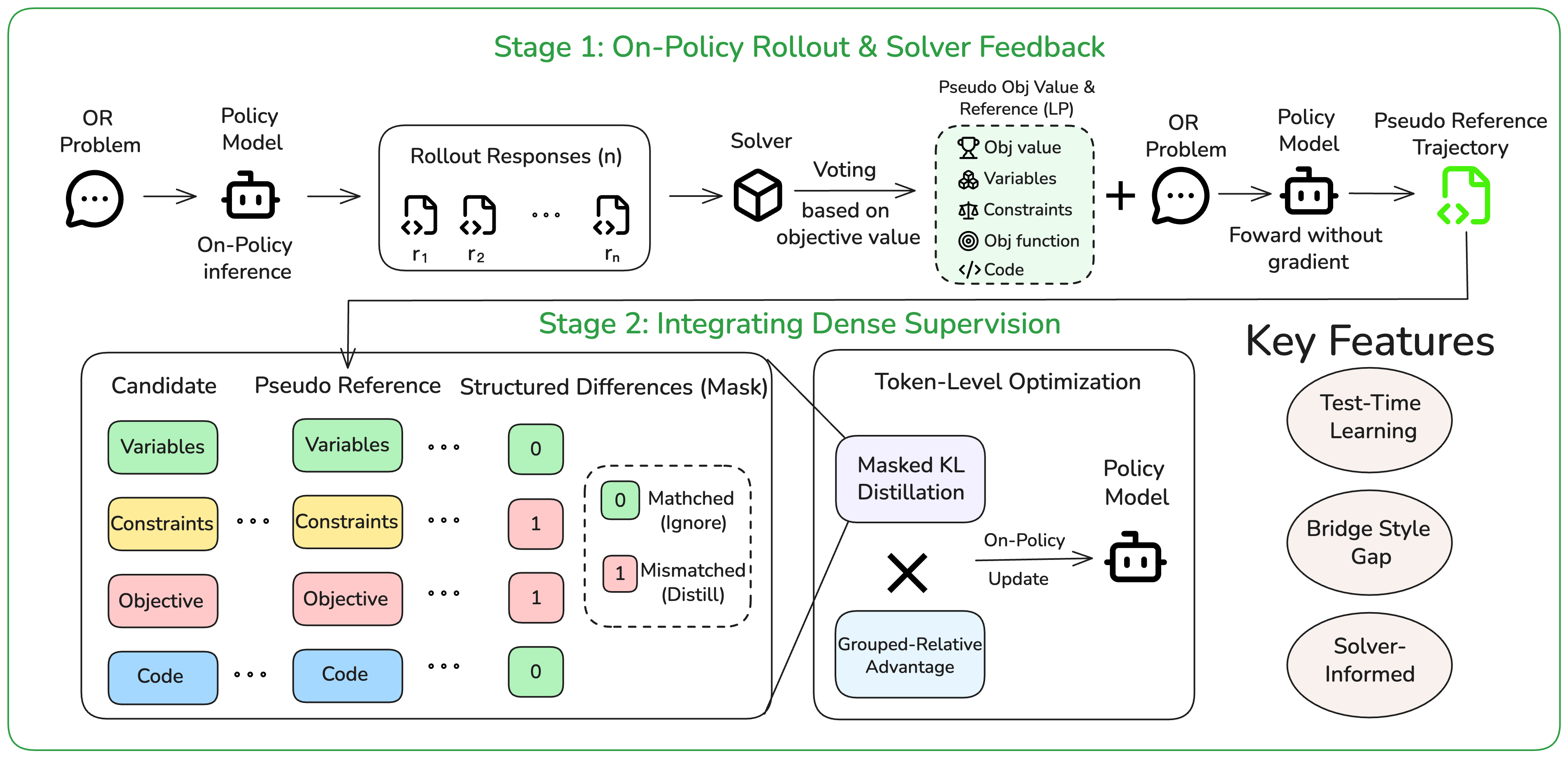}
    \caption{
    Overview of \method{}.
    The policy first generates on-policy trajectories whose programs are
    evaluated and grouped by solver feedback, yielding group-relative rewards
    and a majority-group solver artifact.
    The same policy then uses this artifact as privileged context, while an
    LP-structured mask restricts self-distillation to structurally mismatched
    variable, constraint, objective, and code sections.
    The masked token-level signal and sequence-level reward jointly update the
    deployable policy.
    }
    \label{fig:solid_method}
\end{figure}

\noindent\textbf{Overview.}
\method{} follows the three components summarized in
Figure~\ref{fig:solid_method}. \textbf{On-policy inference} samples $N$
responses under a fixed structured prompt. A \textbf{solver executor and
information extractor} execute the generated code and turn their objective
values and LP artifacts into a majority-LP pseudo-reference and
candidate--reference LP-structure differences. \textbf{Solver-informed
learning} uses solver-outcome agreement for sequence-level advantage and the
reference artifact for self-distillation.

\subsection{\method{}}

\paragraph{On-policy inference.}
A fixed response template makes each rollout structurally separable: it asks the model to reason step by step in ordered sections with declared modeling roles. The system message supplies the format instruction, while the user message contains the OR problem. The reference field is empty during rollout and contains the majority-LP pseudo-reference only during contextual re-scoring. This design is not tied to a particular prompt template: for any separable response schema, solver-grounded sections can be associated with the corresponding artifact evidence and masked according to candidate--reference differences. In our implementation, we instantiate this principle with a nine-stage response template; the complete prompt is provided in supplementary material. An abridged template is:
\begin{center}
\fbox{\begin{minipage}{0.92\columnwidth}
\footnotesize\ttfamily
SYSTEM:\\
\{FORMAT INSTRUCTION\}\\
Reason step by step in separated sections (<step>...</step>).\\
Use the solver-information reference only in allowed sections.\\
\{REFERENCE LP\}\\
(contextual scoring only; empty during rollout)\\[3pt]
USER:\\
\{OR PROBLEM\}
\end{minipage}}
\end{center}
The prompt mapping identifies where solver evidence is applicable; the binary
LP-structured mask itself is determined solely by candidate--reference
artifact comparison.

\paragraph{Solver executor and information extractor.}
For each rollout, the executor records execution status, objective
$\hat z_i$, and LP artifact $a_i$. After normalizing objective direction,
executable objectives within $10^{-6}$ form a vote cluster. Let $C^\star$ be
the largest cluster and $\Rvote$ its median. The executable
rollout closest to $\Rvote$ supplies the majority-LP pseudo-reference
$s^\star$. If none executes,
$\Rvote=\bot$ and no reference is constructed.
The vote reward is
\[
r_i^{\mathrm{vote}}=\mathbf{1}[i\in C^\star].
\]

The extractor next compares each candidate artifact with $s^\star$. It
canonicalizes variable identities and LP record order, then builds signatures for variables
(types, bounds, and constraint-matrix columns), the objective (sense and
coefficients), constraints (sorted rows of coefficients, senses, and
right-hand sides), and the complete normalized LP emitted by the code.
Table~\ref{tab:lp-prompt-map} illustrates how these signatures map to the StepORLM
response template \citep{zhou2026steporlm}; other implementations may adapt
the mapping to their response schema and available solver information.
For the mapped sections $\mathcal K=\{3,4,5,9\}$,
\[
d_{i,k}=\mathbf{1}\!\left[F_k(a_i)\ne F_k(s^\star)\right],
\qquad k\in\mathcal K,
\]
where $F_k$ extracts the corresponding signature. If $a_i$ is missing,
invalid, or unseparable, only $d_{i,9}$ is set to one because only the code
failure is observable. The extracted evidence is therefore
\[
E_i=(\Rvote,\ s^\star,\ d_i).
\]

Let $\kappa_i(t)$ denote the response section containing token $t$. The sample
gate and final LP-structured mask are
\[
\begin{aligned}
u_i&=\mathbf{1}[\Rvote\neq\bot]\,\mathbf{1}[i\notin C^\star],\\
b_{i,t}
&=
\begin{cases}
u_i d_{i,\kappa_i(t)}, & \kappa_i(t)\in\mathcal K,\\
0, & \text{otherwise}.
\end{cases}
\end{aligned}
\]
Thus, majority-cluster responses receive GRPO only, while non-majority
responses can receive self-distillation only in LP-mismatched mapped
sections.

\begin{table}[t]
\centering
\tblfont
\setlength{\tabcolsep}{4pt}
\caption{Mapping LP signatures to the StepORLM response template.}
\label{tab:lp-prompt-map}
\begin{tabular}{p{0.48\columnwidth}p{0.42\columnwidth}}
\toprule
Compared LP signature & Response section \\
\midrule
Variables, types, bounds, or columns
& 3. Decision Variables \\
Objective sense or coefficients
& 4. Objective Function \\
Row coefficients, senses, or RHS
& 5. Constraints \\
Complete normalized LP emitted by the candidate code
& 9. Python Code Using \texttt{gurobipy} \\
\bottomrule
\end{tabular}
\end{table}

\paragraph{Solver-informed learning objective.}
The majority objective and majority-LP pseudo-reference are inserted into a
context-augmented prompt,
\[
\widetilde x_i=g_{\mathrm{OR}}(x,\Rvote,s^\star).
\]
The rollout policy re-scores the same response without gradients under
$\widetilde x_i$, while the updated student sees only the original problem:
\[
\begin{aligned}
\ell^T_{i,t}
&=\log\pi_{\mathrm{old}}
(y_{i,t}\mid\widetilde x_i,y_{i,<t}),\\
\ell^\theta_{i,t}
&=\log\pi_\theta(y_{i,t}\mid x,y_{i,<t}).
\end{aligned}
\]
This same-policy comparison provides token-level credit without verified
answers or another evaluator. Because privileged LP context can also change
notation, ordering, or derivation detail, uniform self-distillation suffers
from \textbf{style mismatch}. The sample gate avoids KL across
consensus-equivalent responses, and the LP-structure gate prevents contextual
changes in matched or unrelated sections.

Let $\Delta_{i,t}=\ell^T_{i,t}-\ell^\theta_{i,t}$. We use the non-negative
$k_3$ estimator \citep{schulman2020approximatingkl},
\[
\phi_{k_3}(\Delta)=e^\Delta-\Delta-1,
\qquad
\mathcal L_{\mathrm{KL}}
=
\mathbb E_{i,t}\!\left[
b_{i,t}\phi_{k_3}(\Delta_{i,t})
\right].
\]
The rollout importance ratio appears only in the GRPO/PPO policy loss and does
not multiply $\mathcal L_{\mathrm{KL}}$. If no majority-LP pseudo-reference
exists, then $b_{i,t}=0$ and the KL term vanishes. The complete learning
objective is
\[
\boxed{\;
\mathcal L_{\mathrm{SOLID}}
=\mathcal L_{\mathrm{GRPO}}+\beta\mathcal L_{\mathrm{KL}}
\;}
\]
GRPO supplies sequence-level competition, while the masked $k_3$ term provides
LP-structured token-level credit. The following section analyzes how
the group advantage anchors these contextual corrections.

\begin{algorithm}[t]
\caption{\method{} on unlabeled OR prompts}
\label{alg:solid}
\begin{algorithmic}[1]
\REQUIRE Prompts $\mathcal D$, policy $\pi_\theta$, rollouts $N$, coefficient $\beta$
\FOR{each minibatch $\mathcal B\subset\mathcal D$}
  \STATE Set rollout snapshot $\pi_{\mathrm{old}}\leftarrow\pi_\theta$.
  \FOR{each $x\in\mathcal B$}
    \STATE \textbf{Inference:} sample $\{y_i\}_{i=1}^{N}\sim\pi_{\mathrm{old}}(\cdot\mid x)$.
    \STATE \textbf{Executor:} extract and run code; record status, $\hat z_i$, and $a_i$.
    \STATE \textbf{Extractor:} form $C^\star$, $\Rvote$, $s^\star$, rewards, and $A_i$.
    \FOR{each rollout $i$}
      \IF{$\Rvote\neq\bot$ and $i\notin C^\star$}
        \STATE Compare $F_k(a_i)$ with $F_k(s^\star)$; form $b_{i,t}$.
        \STATE Cache $\ell^T_{i,t}$ by scoring $y_i$ under $\widetilde x_i$.
      \ELSE
        \STATE Set $b_{i,t}=0$ for all $t$.
      \ENDIF
    \ENDFOR
  \ENDFOR
  \STATE Compute $\mathcal L_{\mathrm{GRPO}}$ and $\mathcal L_{\mathrm{KL}}$.
  \STATE Update $\theta$ using $\mathcal L_{\mathrm{GRPO}}+\beta\mathcal L_{\mathrm{KL}}$.
\ENDFOR
\end{algorithmic}
\end{algorithm}

%% file: sections/theory.tex
\section{Anchored Correction under Style Mismatch}
\label{sec:style-theory}

Contextual likelihood shifts conflate task correction with \textbf{style mismatch} in
wording, order, or derivation detail. \method{} separates their roles: the LP-structured
mask determines where context may act, while the group advantage determines
how strong that preference must be to reverse an update.
Proofs are provided in the supplementary material.
Let $p_{i,t}=\pi_\theta(y_{i,t}\mid x,y_{i,<t})$,
$q_{i,t}=\pi_{\mathrm{old}}(y_{i,t}\mid\widetilde x_i,y_{i,<t})$, and
$\Delta_{i,t}=\log q_{i,t}-\log p_{i,t}$. 
Here we consider the unclipped scenarios.

\noindent\textbf{Lemma 1 (masked contextual direction).}
For fixed $q_{i,t}$ and importance ratios,
\[
-\nabla_\theta\mathcal L_{\mathrm{KL}}
=\mathbb E_{i,t}\!\left[
\widetilde b_{i,t}
\left(\frac{q_{i,t}}{p_{i,t}}-1\right)
\nabla_\theta\log p_{i,t}\right].
\]
% \emph{Proof.}
% Use $\phi'(\Delta)=e^\Delta-1$ and
% $\nabla_\theta\Delta=-\nabla_\theta\log p_{i,t}$.

% At $\theta=\theta_{\mathrm{old}}$, $\widetilde\rho_{i,t}=1$. Inside the
% unclipped GRPO region,
\[
\begin{aligned}
A^{\mathrm{eff}}_{i,t}
&=A_i+\beta b_{i,t}(e^{\Delta_{i,t}}-1),\\
-\nabla_\theta\mathcal L_{\mathrm{SOLID}}
&=\mathbb E_{i,t}
[A^{\mathrm{eff}}_{i,t}\nabla_\theta\log p_{i,t}].
\end{aligned}
\]
For small $\Delta_{i,t}$, this becomes
\[
A^{\mathrm{eff}}_{i,t}
=A_i+\beta b_{i,t}(\ell^T_{i,t}-\ell^S_{i,t})
+O(\beta\Delta_{i,t}^2),
\]
which formalizes the log-probability-difference approximation.

Negative sample reinforcement suppresses a sampled response and redistributes
mass toward alternatives \citep{zhu2025negative}. Non-majority rollouts are
similarly \emph{vote-negative}, though not necessarily incorrect. When
$A_i<0$, the negative samples can therefore suppress a locally useful token.

\noindent\textbf{Proposition 1 (anchored sign reversal).}
For $A_i<0$, $b_{i,t}=0$ gives $A^{\mathrm{eff}}_{i,t}=A_i$. If $b_{i,t}=1$,
let
\[
\tau_i=\log\!\left(1-\frac{A_i}{\beta}\right)>0.
\]
Then $A^{\mathrm{eff}}_{i,t}\leq A_i$ for $\Delta_{i,t}\leq0$;
$A^{\mathrm{eff}}_{i,t}\leq0$ for
$0<\Delta_{i,t}\leq\tau_i$; and $A^{\mathrm{eff}}_{i,t}>0$ only for
$\Delta_{i,t}>\tau_i$.
% \emph{Proof.}
% Solve the sign of $A_i+\beta(e^{\Delta_{i,t}}-1)$.

Thus masking removes the KL channel for majority, matched, and unmapped tokens,
regardless of their surface-form shift. Within an activated section, $A_i$ is
an anchor: weak contextual preferences cannot reverse a negative update, while
a strong LP-informed preference can preserve a locally useful token. 
The LP-structured mask ensures that this occurs only at error-related reasoning stages, thereby preventing reversals caused by \textbf{style mismatch}.

%% file: sections/experiments.tex
\section{Experiments}

\subsection{Experimental Setup}

\paragraph{Benchmarks and Metrics.}%TODO:说这些是challenging的benchmark @cml
We evaluate on OptMATH \citep{lu2025optmath}, MAMO-Complex \citep{huang2024mamo}, and IndustryOR (InOR) \citep{huang2025orlm}, three particularly demanding benchmarks spanning general mathematical optimization, complex executable modeling, and real-world industrial OR scenarios. These benchmarks require end-to-end semantic, mathematical, and executable correctness rather than surface-form matching or short-answer accuracy.
We report majority accuracy maj@$N$ and pass@$k$. Maj@$N$ denotes the accuracy of the answer selected by majority voting among $N$ sampled solutions. Following \citet{chen2021codex}, pass@$k$ measures the probability that at least one of $k$ sampled solutions is correct and is estimated as
\[
\operatorname{pass@}k =
1-\frac{\binom{n-c}{k}}{\binom{n}{k}},
\]
where $n$ is the total number of generated samples and $c$ is the number of correct samples. We mark a sample as correct if its derived objective value is within $10^{-6}$ of the ground-truth value. All values reported in the tables are percentages. 

\paragraph{Baselines.}
\begin{itemize}
\item \textbf{Base.} We evaluate both Qwen3-4B-Instruct \citep{yang2025qwen3} and StepORLM \citep{zhou2026steporlm}, an OR-specialized model built on Qwen3-8B. The corresponding model before task-specific training is reported as \textsc{Base}.

\item \textbf{TTRL.} The outcome-only, group-relative baseline that derives rewards from majority voting over objective values\citep{zuo2025ttrl}.
\item \textbf{OPSD.} A whole-response self-distillation baseline in which
the policy re-scores its own rollouts under privileged LP context, without
LP-structured error localization \citep{zhao2026selfdistilled}. Unless
otherwise stated, OPSD uses the majority-LP solution as the
privileged information.
\item \textbf{\method{}.} Our method combines the TTRL-style
group-relative objective with solver-informed, LP-structured masked
self-distillation.

\end{itemize}

\paragraph{Implementation details.}
For the main training run, we construct a 50,000-instance subset from the publicly released \textbf{OptMATH-Train} corpus \citep{lu2025optmath,auroragem2025optmathtrain}. Specifically, we generate five independent responses to each candidate problem using Qwen3-8B and retain problems that are solved correctly at least once but no more than three times, i.e., with a success count in $\{1,2,3\}$ out of five trials.

Unless otherwise specified, all generated programs in our pipelines are executed using Gurobi
\citep{gurobi2026}. Gurobi provides the solver statuses and artifacts used to determine the LP-structured mask. For StepORLM \citep{zhou2026steporlm}, we retain its native COPT backend \citep{ge2023copt} to avoid introducing conversion errors between solver APIs.

During training, we sample 24 trajectories per prompt at a temperature of 1.0 with a batch size of 32. We use a GRPO clipping range of 0.2, set the KL coefficient to $10^{-3}$, and allow responses of up to 32,768 tokens. Training is conducted on four NVIDIA H100 GPUs ($4\times$H100), and our implementation is built on \textsc{veRL} \citep{sheng2024hybridflow}. Evaluation keeps the same configuration as well and checkpoints are evaluated at the same training steps.

\subsection{Main Results}

\begin{table}[tb]
\centering
\caption{Base and matched-budget results on three datasets. TTRL and
\method{} use the same training budget and evaluation protocol. Bold marks
the best value within each dataset--metric pair.}
\label{tab:paired-main-results}
\begin{subtable}[t]{0.492\textwidth}
\centering
\caption{Qwen3-4B-Instruct.}
\label{tab:qwen-main-results}
\scriptsize
\setlength{\tabcolsep}{1.5pt}
\resizebox{\linewidth}{!}{%
\begin{tabular}{@{}llrrrr@{}}
\toprule
Dataset & Method & maj@64 & pass@1 & pass@2 & pass@4 \\
\midrule
OptMATH & Base & 24.70 & 14.16 & 21.89 & 31.33 \\
OptMATH & TTRL & 30.12 & 19.01 & 25.29 & 31.32 \\
OptMATH & \method{} & \textbf{39.76} & \textbf{24.25} &
\textbf{31.62} & \textbf{38.70} \\
\midrule
MAMO-Complex & Base & \textbf{35.47} & 25.00 & 34.81 & \textbf{44.33} \\
MAMO-Complex & TTRL & 33.50 & 27.56 & 34.95 & 41.52 \\
MAMO-Complex & \method{} & 33.00 & \textbf{29.51} &
\textbf{35.71} & 41.34 \\
\midrule
InOR & Base & 52.00 & 41.56 & 48.84 & 55.11 \\
InOR & TTRL & \textbf{53.00} & 43.84 & 50.40 & 55.95 \\
InOR & \method{} & \textbf{53.00} & \textbf{43.92} &
\textbf{50.85} & \textbf{56.64} \\
\bottomrule
\end{tabular}%
}
%\caption{Qwen3-4B-Instruct.}
%\label{tab:qwen-main-results}
\end{subtable}\hfill
\begin{subtable}[t]{0.492\textwidth}
\centering
\caption{StepORLM.}
\label{tab:steporlm-main-results}
\scriptsize
\setlength{\tabcolsep}{1.5pt}
\resizebox{\linewidth}{!}{%
\begin{tabular}{@{}llrrrr@{}}
\toprule
Dataset & Method & maj@64 & pass@1 & pass@2 & pass@4 \\
\midrule
OptMATH & Base & 17.47 & 10.99 & 16.67 & 22.29 \\
OptMATH & TTRL & 30.72 & 17.48 & 23.64 & 30.38 \\
OptMATH & \method{} & \textbf{31.33} & \textbf{18.25} &
\textbf{24.40} & \textbf{30.53} \\
\midrule
MAMO-Complex & Base & 62.56 & 52.22 & 65.02 & 72.91 \\
MAMO-Complex & TTRL & 68.97 & 61.22 & 68.08 & 72.26 \\
MAMO-Complex & \method{} & \textbf{70.44} & \textbf{66.43} &
\textbf{71.58} & \textbf{74.79} \\
\midrule
InOR & Base & 46.00 & 36.50 & 43.90 & 49.25 \\
InOR & TTRL & \textbf{48.00} & 39.36 & 46.19 & \textbf{51.21} \\
InOR & \method{} & \textbf{48.00} & \textbf{39.81} &
\textbf{47.03} & 50.59 \\
\bottomrule
\end{tabular}%
}
%\caption{StepORLM.}
%\label{tab:steporlm-main-results}
\end{subtable}
%\caption{Base and matched-budget results on three datasets. TTRL and
%\method{} use the same training budget and evaluation protocol. Bold marks
%the best value within each dataset--metric pair.}
%\label{tab:paired-main-results}
\end{table}

\paragraph{Qwen3-4B-Instruct results.}
Table~\ref{tab:qwen-main-results} compares Base, TTRL, and \method{} under
the same training budget. The strongest gains occur on OptMATH, where majority accuracy rises by 9.64 points and pass@$k$ improves by an average of 6.31 points across $k\in\{1,2,4\}$. Improvements also extend across all three pass@$k$ metrics on InOR. When OptMATH and InOR are aggregated over their 266 problems, \method{} improves majority accuracy by 6.02 points and pass@$k$ by an average of 4.09 points over TTRL. The results reveal steady improvements of \method{} across datasets.

\paragraph{StepORLM results.}
Table~\ref{tab:steporlm-main-results} compares Base, TTRL, and \method{} on an OR-tuned model. Under the same training budget, \method{} improves or matches majority accuracy on all three datasets. It also outperforms TTRL on eight of twelve metrics overall. The largest gains occur on MAMO-Complex, where pass@$k$ improves by an average of 3.75 points across $k\in{1,2,4}$. On OptMATH and InOR, the pass@$k$ results remain broadly comparable to TTRL, with gains at smaller $k$ and slight decreases at larger $k$.

\subsection{Reasoning Capacity of the Prompted Student}

\begin{table}[t]
\centering
\tblfont
\setlength{\tabcolsep}{5.0pt}
\caption{Reference-source diagnostic for Qwen3-4B-Instruct on InOR using
64 samples per problem. Correct references use
evaluation labels only to measure correction capacity. Arrows indicate
changes relative to ``None''; bold marks the best value in each column.}
\label{tab:reference-ablation}
\begin{tabular}{@{}lrrrr@{}}
\toprule
Reference & maj@64 & pass@1 & pass@2 & pass@4 \\
\midrule
None
& 52.00 & 41.56 & 48.84 & 55.11 \\
Correct
& \textbf{63.00}$\uparrow$ & \textbf{49.59}$\uparrow$
& \textbf{55.50}$\uparrow$ & \textbf{59.82}$\uparrow$ \\
Majority
& 51.00$\downarrow$ & 43.11$\uparrow$
& 47.37$\downarrow$ & 50.98$\downarrow$ \\
Random
& 52.00 & 40.50$\downarrow$
& 47.27$\downarrow$ & 52.76$\downarrow$ \\
\bottomrule
\end{tabular}
\end{table}

Table~\ref{tab:reference-ablation} tests which information source is useful in
the privileged-context prompt. The correct reference provides an upper bound on the model’s capacity to exploit reliable solver information and improves its response accuracy. 
It can be inferred that the prompted
student can exploit sound solver information, but correction is sensitive to
reference quality. This gap motivates combining the group-relative objective with localized self-distillation. The group-relative objective provides a stable learning signal that strengthens the model, enabling it to generate more reliable supervision for subsequent self-distillation.

\subsection{Ablation with Different Masking Strategies}
\paragraph{Masking strategies.}
Table~\ref{tab:qwen-kl-ablation} gives a matched-budget comparison among
different KL masking strategies. Overall, \method{} provides the strongest
cross-dataset balance, leading most metrics on OptMATH and MAMO-Complex while
remaining competitive on InOR. Its broad gains across majority accuracy and
pass@$k$ suggest that LP-structured localization yields more reliable
self-distillation than applying KL uniformly or through a random section
mask. TTRL has no KL term, TTRL + KL applies whole-response KL, and random KL
masks every response section independently with probability 0.5.

\begin{table}[H]
\centering
\footnotesize
\setlength{\tabcolsep}{5.5pt}
\renewcommand{\arraystretch}{1.02}
\caption{Qwen3-4B-Instruct masking ablation under a matched training budget.
Bold marks the best and underline marks the second-best result within each
dataset.}
\label{tab:qwen-kl-ablation}
\begin{tabular*}{0.84\textwidth}{@{\extracolsep{\fill}}llrrrr@{}}
\toprule
Dataset & Method & maj@64 & pass@1 & pass@2 & pass@4 \\
\midrule
\multirow{4}{*}{OptMATH}
 & TTRL & 30.12 & 19.01 & 25.29 & 31.32 \\
 & TTRL + KL & \textbf{39.76} & \underline{21.20} & \underline{28.81} & \underline{36.61} \\
 & TTRL + random KL & 30.12 & 18.52 & 24.45 & 29.83 \\
 & \method{} & \textbf{39.76} & \textbf{24.25} & \textbf{31.62} & \textbf{38.70} \\
\midrule
\multirow{4}{*}{MAMO-Complex}
 & TTRL & 33.50 & 27.56 & 34.95 & 41.52 \\
 & TTRL + KL & \underline{33.99} & \textbf{31.47} & \underline{36.70} & \underline{41.70} \\
 & TTRL + random KL & 33.00 & 30.70 & 36.01 & 40.97 \\
 & \method{} & \textbf{35.96} & \underline{31.02} & \textbf{36.76} & \textbf{42.22} \\
\midrule
\multirow{4}{*}{InOR}
 & TTRL & 53.00 & 43.84 & 50.40 & 55.95 \\
 & TTRL + KL & \underline{54.00} & 43.42 & 50.04 & 55.62 \\
 & TTRL + random KL & \textbf{55.00} & \textbf{46.28} & \textbf{52.54} & \textbf{57.54} \\
 & \method{} & 53.00 & \underline{45.13} & \underline{51.47} & \underline{56.64} \\
\bottomrule
\end{tabular*}
\end{table}

\begin{figure}[H]
\centering
\begin{subfigure}[t]{0.48\linewidth}
  \centering
  \includegraphics[width=\linewidth]{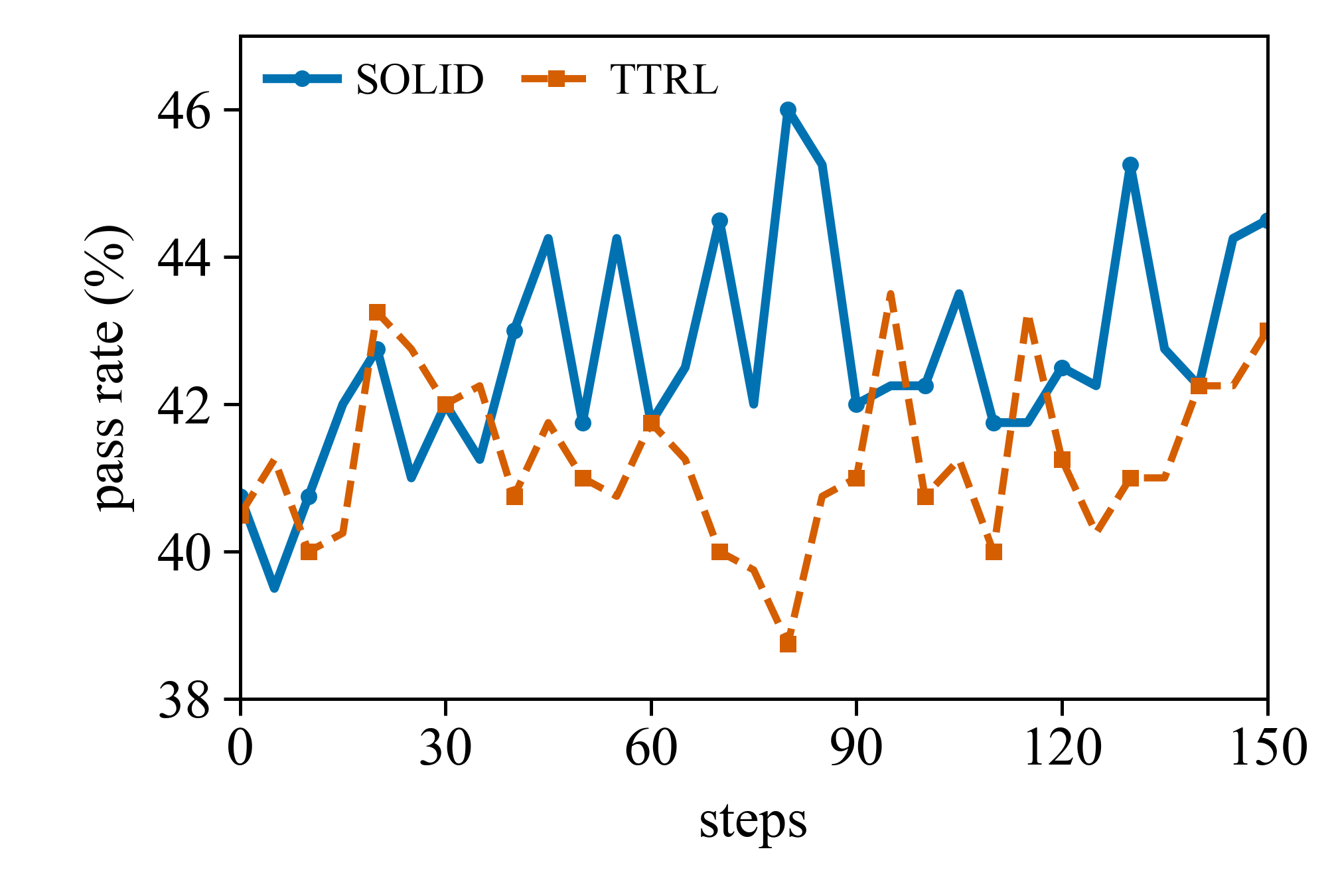}
  \caption{\method{} versus outcome-only TTRL.}
  \label{fig:solid-vs-ttrl-training}
\end{subfigure}\hfill
\begin{subfigure}[t]{0.48\linewidth}
  \centering
  \includegraphics[width=\linewidth]{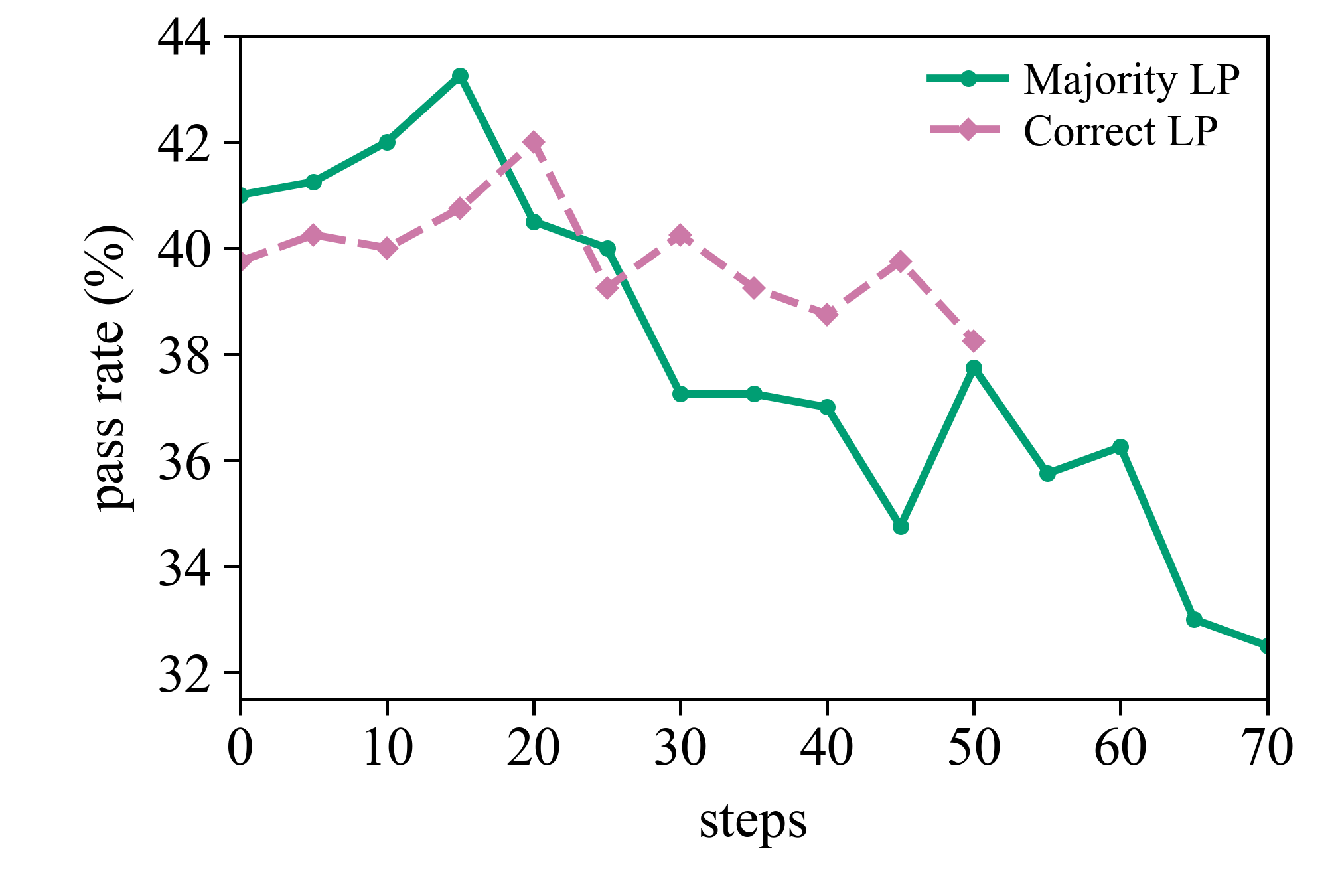}
  \caption{Whole-response OPSD with two LP contexts.}
  \label{fig:opsd-failure}
\end{subfigure}
\caption{Training dynamics of Qwen3-4B-Instruct: validation reports mean@4 performance on IndustryOR.
(a) The complete \method{} objective maintains an advantage over
outcome-only TTRL. (b) Pure whole-response OPSD declines after an early peak
with both majority-LP and diagnostic correct-LP context.}
\label{fig:training-dynamics}
\end{figure}

\paragraph{Training dynamics and OPSD failure.}
Figure~\ref{fig:solid-vs-ttrl-training} compares \method{} directly with
TTRL and shows a clearer advantage across steps for the \method{} objective.
Figure~\ref{fig:opsd-failure} isolates the failure of pure OPSD:
mean@4 first rises and then declines under both majority-LP and correct-LP
privileged context, with a sharper decline for the majority-LP run. Because
the correct-LP curve uses evaluation labels only as a diagnostic, its decline
shows that whole-response self-distillation can fail even after reference
error is removed. This is consistent with privileged-context style mismatch,
although the trajectory does not establish style mismatch as the unique
cause.

\paragraph{Further discussion.}
Several extensions could further strengthen SOLID. First, more reliable extraction of candidate LPs could reduce erroneous section activation when solver artifacts are incomplete or unstable. Second, a more complicated aggregation algorithm could make the majority-LP pseudo-reference more reliable when rollout consensus is incorrect. Finally, replacing the binary LP-structured mask with an adaptive, confidence-weighted localization mechanism may better preserve useful reasoning while limiting context-induced style shifts, particularly on datasets and metrics where the current gains are less consistent.

%% file: sections/conclusion.tex
\section{Conclusion}

We presented \method{}, a post-training framework that bootstraps LLMs from
unlabeled OR problems using on-policy rollouts and solver outcomes. \method{}
addresses three central obstacles to modern OR model training: costly scalable
supervision, fine-grained credit assignment and style mismatch. It combines solver-outcome
sequence-level learning with same-policy token-level self-distillation, using
the binary LP-structured mask to localize guidance to the mismatched reasoning parts. Improvements across both general-purpose and
OR-tuned models show that this signal remains useful before and after
domain-specific tuning. Furthermore, our OPSD diagnostic exposes a failure mode of
self-distillation in OR and motivates
\method{}'s LP-structure masking. Thus, these designs lay a solid foundation for future work on self-distillation for OR-oriented LLM training.

%% file: sections/appendix.tex
% !TeX root = appendix_cml_standalone.tex
% Technical appendix for SOLID.
% AAAI submission note: submit the compiled appendix PDF as supplementary
% material, not as part of the main paper PDF. The main paper should remain
% self-contained and may cite this appendix only for supporting details.
% Items marked TODO(CML) require experiment logs or an author decision and are
% intentionally hidden from the compiled PDF.

\appendix
\setcounter{secnumdepth}{2}
\renewcommand{\thesection}{Appendix \Alph{section}}
\renewcommand{\thesubsection}{\Alph{section}.\arabic{subsection}}

\section{Gurobi Prompts}
\label{app:prompts}
\subsection{Training and Inference Prompt}
\label{app:training-inference-prompt}

The same Gurobi system and user messages are used for training rollouts and
inference.

\begingroup
\tiny
\noindent\hrule\smallskip
\setlength{\columnsep}{18pt}
\setlength{\columnseprule}{0.2pt}
\begin{multicols}{2}
\raggedcolumns
\begin{verbatim}
SYSTEM
You are a highly specialized AI assistant with deep expertise
in mathematical modeling, Python programming, and the Gurobi
solver. Your primary mission is to transform user-provided
optimization problems into clear, structured, and solvable
models.

When a user presents an optimization question, you must
rigorously analyze it and deliver a comprehensive response. To
ensure maximum clarity, consistency, and correctness, your
entire output must strictly adhere to the following nine-step
structure. Do not add, omit, or reorder these steps.

**Your Response Structure:**

1. **Problem Description**: Concisely summarize the user's
problem in your own words.
2. **Sets and Parameters**: Define all the sets, indices, and
known parameters.
3. **Decision Variables**: Clearly define the variables the
model will solve for.
4. **Objective Function**: State the objective function with a
clear explanation of its purpose.
5. **Constraints**: Detail each constraint with a brief
explanation of what it represents.
6. **Mathematical Model**: Present the complete mathematical
formulation using clear notation.
7. **Nonlinear Relationships**: If any, describe nonlinearities
and how they will be handled (e.g., linearization). If none,
state "The model is linear."
8. **Final Model**: Present the final, complete mathematical
model ready for implementation.
9. **Python Code Using `gurobipy`**: Provide a complete and
executable Python script that uses the `gurobipy` library to
solve the model. The code should be well-commented to link back
to the mathematical formulation.

**Formatting Instructions:**

* You **must** output exactly nine `<step>...</step>` blocks, in
the order listed above.
* Each block must start with the exact step title in bold as the
first line inside the tag.
* Put the analysis or implementation content after the title,
inside the same `<step>...</step>` block.
* Use the standard closing tag `</step>`. Do not use `<\step>`,
`</ Step>`, or any other variant.
* Do not wrap the `<step>...</step>` blocks in bullets, numbered
lists, or any outer container.
* Use LaTeX formatting for mathematical notation, enclosing
formulas in `$` or `$$` delimiters.
* In the `Python Code Using gurobipy` step, put the complete
executable code inside a ```python fenced code block.

**Required output skeleton:**

<step>
**Problem Description**
Summarize the optimization problem.
</step>

<step>
**Sets and Parameters**
Define all sets, indices, and parameters.
</step>

<step>
**Decision Variables**
Define every decision variable and its domain.
</step>

<step>
**Objective Function**
State the objective and whether it is minimized or maximized.
</step>

<step>
**Constraints**
List and explain all constraints.
</step>

<step>
**Mathematical Model**
Give the complete formulation.
</step>

<step>
**Nonlinear Relationships**
State whether the model is linear; if not, explain the
linearization.
</step>

<step>
**Final Model**
Present the final model ready for implementation.
</step>

<step>
**Python Code Using gurobipy**
```python
import gurobipy as gp
from gurobipy import GRB

model = gp.Model("model")
# build variables, objective, and constraints
model.optimize()

if model.status == GRB.OPTIMAL:
    solution = {var.VarName: var.X for var in model.getVars()}
    print("Just print the best obj:", model.ObjVal)
else:
    print("No Solution")
```
</step>

**Gurobi Code Requirements:**
* Make sure to import necessary packages, such as
"import gurobipy as gp" and "from gurobipy import GRB".
* When you create a model, make sure to use
"model = gp.Model("model")".
* When you add a variable, use "vtype=GRB.CONTINUOUS",
"vtype=GRB.INTEGER", or "vtype=GRB.BINARY".
* Do not name variables and constraints.
* Use "model.addConstr()" or "model.addConstrs()" to add
constraints.
* If you want to set "lb" or "ub" as infinity, please use
"lb=-GRB.INFINITY" or "ub=GRB.INFINITY".
* When you set objective, you should use the
"model.setObjective" method and use "GRB.MINIMIZE" or
"GRB.MAXIMIZE".
* Make sure to use "model.optimize()" to solve the question.
* The code output statement is:
if model.status == GRB.OPTIMAL:
    solution = {{var.VarName: var.X for var in model.getVars()}}
    print("Just print the best obj:", model.ObjVal)
else:
    print("No Solution")

Begin your work once the user provides the optimization problem.

USER
Below is an optimization modeling question. Build a mathematical
model and corresponding python code using gurobipy that
appropriately addresses the question:

{question}

Think step by step.
\end{verbatim}
\end{multicols}
\smallskip\hrule
\endgroup

\section{Benchmark Details}
\label{app:benchmarks}

We evaluate on three complementary benchmarks. OptMATH provides the
original 166-instance benchmark for broad mathematical-optimization modeling,
whereas we evaluate MAMO-Complex and IndustryOR using the corrected releases
provided by SIRL \citep{NEURIPS2025_98555b92}. Table~\ref{tab:appendix-benchmarks}
summarizes the exact evaluation sets used in our experiments.

\begin{center}
\centering
\captionof{table}{Evaluation benchmark sizes.}
\label{tab:appendix-benchmarks}
\small
\setlength{\tabcolsep}{5pt}
\begin{tabular}{lc}
\toprule
Benchmark & Instances \\
\midrule
OptMATH & 166 \\
MAMO-Complex (SIRL-corrected) & 203 \\
IndustryOR (SIRL-corrected) & 100 \\
\bottomrule
\end{tabular}
\end{center}

% TODO(CML): Add the exact SIRL repository commit or data-release date used for
% MAMO-Complex and IndustryOR.

\subsection{OptMATH}
\label{app:benchmark-optmath}

The OptMATH benchmark was introduced by Lu et al. as a scalable benchmark for
optimization modeling built through bidirectional data synthesis and rejection
filtering \citep{lu2025optmath}. Its benchmark split contains long natural
language problem descriptions and spans multiple mathematical program classes,
including LP, MILP, IP, NLP, SOCP, and related optimization models. In our
experiments, we use the original 166-instance OptMATH benchmark without
additional filtering or correction. We draw the training pool separately from
OptMATH-Train and summarize run-level details in the reproducibility notes
below.

% TODO(CML): State the exact split/version or commit and any excluded instances.

\subsection{MAMO-Complex}
\label{app:benchmark-mamo}

The MAMO dataset was introduced by Huang et al. as a benchmark for evaluating
large language models on mathematical modeling with solver-executable answers
\citep{huang2024mamo}. It contains two main LP-oriented components, EasyLP and
ComplexLP, with ComplexLP focusing on more involved formulations and longer
modeling chains. The original ComplexLP subset contains 211 instances. Following
SIRL, we use the corrected MAMO-ComplexLP release, where invalid or ambiguous
instances were removed or revised, resulting in 203 evaluated problems
\citep{NEURIPS2025_98555b92}. We use this corrected split in full for our
MAMO-Complex experiments.

% TODO(CML): If available, add the filename used by the evaluation script, e.g.,
% MAMO_ComplexLP_fixed.jsonl.

\subsection{IndustryOR}
\label{app:benchmark-inor}

The IndustryOR dataset, introduced with ORLM, is an industrial operations
research benchmark for testing LLMs on practical domain-specific optimization
tasks \citep{huang2025orlm}. It contains 100 real-world scenarios drawn from 13
industries and covers five OR task categories across three difficulty levels.
Following the SIRL corrected release, we use the revised IndustryOR benchmark
with expert-reviewed questions and answers \citep{NEURIPS2025_98555b92}. The
correction preserves the full 100-instance benchmark, and we include all 100
problems in our IndustryOR evaluation.

% TODO(CML): Add category-level instance counts if they are available.

\subsection{Correctness Criterion}
\label{app:correctness}

For each generated response, we extract and execute the solver program. For
benchmarks with numeric targets, we count a sample as correct only when
execution completes and the extracted objective \(\hat z\) satisfies
\[
    \frac{\left|\hat z-z^\star\right|}
         {\left|z^\star\right|+1} < 10^{-6}.
\]
We score missing code, extraction failures, timeouts, runtime exceptions,
nonnumeric objectives, and non-optimal solver outcomes as incorrect. We use
this ground-truth criterion only for evaluation and construct the training vote
from rollout objectives without using \(z^\star\).

\FloatBarrier

\section{Additional Reproducibility Details}
\label{app:reproducibility}

\subsection{Software and Inference Stack}
\label{app:software-stack}

We vendor a \textsc{veRL} 0.7.0.dev codebase and use asynchronous vLLM rollout
with the vLLM V1 engine enabled. We require Python 3.10 or newer, PyTorch in
\([2.6,2.10)\), vLLM in
\([0.8.5,0.16)\), \texttt{gurobipy} 12.0 or newer, and \texttt{coptpy} 7.2 or
newer. We record these version constraints in the artifact, but do not include
a package-lock file with exact patch versions. We use Gurobi for
Qwen3-4B-Instruct, whereas StepORLM retains COPT.

Beyond the main-paper settings, we use bfloat16 rollout, tensor parallelism of
one, a maximum model length of 38,912 tokens, and a vLLM GPU-memory utilization
target of 0.30. We cap the teacher prompt at 2,048 tokens and right-truncate it
before appending the sampled response. We optimize the actor with FSDP through
the vendored \textsc{veRL} runtime.

\subsection{Code Extraction and Instrumentation}
\label{app:code-extraction}

We apply a deterministic extraction cascade. We first search for code inside
\texttt{<python>...</python>}, then for the first global
\texttt{python}/\texttt{py} Markdown fence, an untyped fence that resembles
Python, or raw text beginning with a solver import. If these paths fail, we
search the titled ninth \texttt{<step>} block and finally the last available
\texttt{<step>} block. We mark a response with no extractable solver code as
failed.

After extraction, we remove Markdown fences. We locate the model variable from
the first \texttt{model.optimize()} or \texttt{model.solve()} call and inject
solver-specific instrumentation. We print the objective only under
\texttt{GRB.OPTIMAL} or \texttt{COPT.OPTIMAL}, print the solution vector used
by diagnostics, and write the final model to an LP file. We recognize the injected
\texttt{Just print the best obj:} line, with two legacy output patterns
retained as fallbacks.

\subsection{Execution Limits and Failure Handling}
\label{app:execution-limits}

We run generated programs in separate worker processes with a 30-second timeout
per program and cap concurrent reward-execution workers at 64. We record
timeouts as \texttt{Timeout Error}. For exceptions, missing code, and process
failures, we produce non-success reports and assign zero execution and answer
credit. We do not set an explicit per-program CPU or memory quota; the
surrounding Ray job and cluster scheduler govern these resources. Accordingly,
our execution layer provides process isolation and timeout control but should
not be interpreted as a hardened security sandbox.

\subsection{Objective Voting and Reference Selection}
\label{app:voting-details}

We admit only executions with status \texttt{Done} and a finite numeric
objective into the vote. We normalize objective direction from the parsed LP
artifact and group values by complete-linkage intervals of absolute width
\(10^{-6}\), preventing tolerance chaining. We rank clusters deterministically
by decreasing size, increasing mean absolute deviation from their median, and
earliest rollout index. Within the winning cluster, we select the usable LP
closest to the median, with rollout order as the final tie-break. If no usable
execution or LP exists, we do not construct a contextual teacher update.

\section{Proofs for Anchored Correction}
\label{app:proofs}

This section supplies the derivations omitted from the main paper and states
the conditions under which its effective-advantage expression is exact. For a
sampled response token, define
\[
\begin{aligned}
p_{i,t}(\theta)
  &=\pi_\theta(y_{i,t}\mid x,y_{i,<t}),\\
p^{\mathrm{old}}_{i,t}
  &=\pi_{\mathrm{old}}(y_{i,t}\mid x,y_{i,<t}),\\
q_{i,t}
  &=\pi_{\mathrm{old}}(y_{i,t}\mid\widetilde x_i,y_{i,<t}),\\
\Delta_{i,t}
  &=\log q_{i,t}-\log p_{i,t}(\theta).
\end{aligned}
\]
We assume positive probabilities on sampled tokens. During an actor update,
the sampled tokens, teacher probabilities $q_{i,t}$, advantages $A_i$, masks,
and averaging weights are held fixed. This is the stop-gradient convention
used in the analysis. Let $b_{i,t}\in[0,1]$ denote the complete fixed gate,
including the non-majority-sample gate and the LP-section gate.

\subsection{Proof of Lemma 1}
\label{app:proof-lemma}

\noindent\textbf{Lemma 1 (Masked contextual direction).}
Under the fixed-teacher and fixed-mask convention above, and away from an
implementation clamp, the masked contextual channel has direction
\[
-\nabla_\theta\mathcal L_{\mathrm{KL}}
=\mathbb E_{i,t}\!\left[
b_{i,t}\left(\frac{q_{i,t}}{p_{i,t}}-1\right)
\nabla_\theta\log p_{i,t}\right].
\]
To be specific, the expectation is taken over on-policy samples, where (i) indexes the sampled trajectory (or response) and (t) indexes the token position within that trajectory. Thus, the pair $(i,t))$ identifies a state-action decision actually visited by the current policy, reflecting the policy’s own behavioral distribution.
With fixed nonnegative averaging weights, the masked contextual loss is
\[
\mathcal L_{\mathrm{KL}}
=\mathbb E_{i,t}\!\left[
b_{i,t}\phi_{k_3}(\Delta_{i,t})\right],
\qquad
\phi_{k_3}(\Delta)=e^\Delta-\Delta-1.
\]
The estimator is nonnegative because $e^u\geq1+u$. Moreover,
$\phi'_{k_3}(\Delta)=e^\Delta-1$, and the fixed teacher gives
\[
\nabla_\theta\Delta_{i,t}
=-\nabla_\theta\log p_{i,t}.
\]
Applying the chain rule gives
\[
\begin{aligned}
\nabla_\theta\mathcal L_{\mathrm{KL}}
&=\mathbb E_{i,t}\left[
b_{i,t}\phi_{k_3}'(\Delta_{i,t})
\nabla_\theta\Delta_{i,t}\right]\\
&=-\mathbb E_{i,t}\left[
b_{i,t}\left(e^{\Delta_{i,t}}-1\right)
\nabla_\theta\log p_{i,t}\right].
\end{aligned}
\]
Since $e^{\Delta_{i,t}}=q_{i,t}/p_{i,t}$,
\[
-\nabla_\theta\mathcal L_{\mathrm{KL}}
=\mathbb E_{i,t}\left[
b_{i,t}\left(\frac{q_{i,t}}{p_{i,t}}-1\right)
\nabla_\theta\log p_{i,t}\right],
\]
which proves the lemma. The contextual term increases the sampled token's
score when $q_{i,t}>p_{i,t}$, decreases it when $q_{i,t}<p_{i,t}$, and
contributes exactly zero when $b_{i,t}=0$.

\paragraph{Reverse-KL interpretation.}
If $y\sim p_\theta$, then
\[
\mathbb E_{y\sim p_\theta}
\!\left[\phi_{k_3}\!\left(\log\frac{q(y)}{p_\theta(y)}\right)\right]
=\mathrm{KL}(p_\theta\|q).
\]
Thus $k_3$ is an unbiased sampled estimator of reverse KL under exact student
sampling. With replayed samples from
$\pi_{\mathrm{old}}\neq\pi_\theta$ and no importance correction on this term,
it remains the implemented nonnegative surrogate but is not generally an
unbiased estimator of that KL divergence. The interpretations coincide at
the on-policy snapshot.

\subsection{Exact Combined Gradient and the Manuscript Form}
\label{app:effective-advantage}

Let
\[
\rho_{i,t}=\frac{p_{i,t}(\theta)}{p^{\mathrm{old}}_{i,t}}
\]
be the GRPO importance ratio. On an unclipped branch,
\[
\mathcal L_{\mathrm{PG}}^{\mathrm{unc}}
=-\mathbb E_{i,t}\!\left[\rho_{i,t}A_i\right].
\]
When the policy-gradient and KL terms use the same token-averaging measure,
their exact combined negative gradient is
\[
-\nabla_\theta
\left(\mathcal L_{\mathrm{PG}}^{\mathrm{unc}}
+\beta\mathcal L_{\mathrm{KL}}\right)
=\mathbb E_{i,t}\!\left[
A^{\mathrm{eff}}_{i,t}\nabla_\theta\log p_{i,t}\right],
\]
where
\begin{equation}
A^{\mathrm{eff}}_{i,t}
=\rho_{i,t}A_i
+\beta b_{i,t}\left(e^{\Delta_{i,t}}-1\right).
\label{eq:app-general-effective}
\end{equation}
If the two token means have denominators $Z_{\mathrm{PG}}$ and
$Z_{\mathrm{KL}}$, respectively, expressing both under the policy-gradient
measure replaces the mask by
\[
\widetilde b_{i,t}
=\frac{Z_{\mathrm{PG}}}{Z_{\mathrm{KL}}}b_{i,t}.
\]
This accounts for the effective mask $\widetilde b_{i,t}$ in Lemma 1 of the
main paper. The unscaled $b_{i,t}$ is valid when both losses share a
normalizer or when the scale is absorbed into $\beta$.

At a strict on-policy actor step, $\rho_{i,t}=1$. Compatible normalizers then
reduce Equation~\eqref{eq:app-general-effective} to
\[
A^{\mathrm{eff}}_{i,t}
=A_i+\beta b_{i,t}\left(e^{\Delta_{i,t}}-1\right).
\]
More generally, this expression is exact if $\rho_{i,t}$ has already been
absorbed into the policy-gradient coefficient denoted by $A_i$.

\paragraph{Small-shift approximation.}
Taylor expansion gives
\[
e^\Delta-1=\Delta+\frac{1}{2}\Delta^2+O(\Delta^3).
\]
Under the manuscript conditions above,
let $\ell^T_{i,t}=\log q_{i,t}$ and
$\ell^S_{i,t}=\log p_{i,t}$. Then
\[
\begin{aligned}
A^{\mathrm{eff}}_{i,t}
&=A_i+\beta b_{i,t}\Delta_{i,t}
+O\!\left(\beta b_{i,t}\Delta_{i,t}^2\right)\\
&=A_i+\beta b_{i,t}
\left(\ell^T_{i,t}-\ell^S_{i,t}\right)
+O\!\left(\beta\Delta_{i,t}^2\right),
\end{aligned}
\]
where the second remainder uses $0\leq b_{i,t}\leq1$. This is a local
approximation and should not be extrapolated to large $|\Delta_{i,t}|$.

\subsection{Proof of Proposition 1}
\label{app:proof-proposition}

\noindent\textbf{Proposition 1 (Anchored sign reversal).}
Under the on-policy manuscript form, let $A_i<0$ and $\beta>0$. If
$b_{i,t}=0$, the contextual channel cannot change the negative coefficient.
If $b_{i,t}=1$, the coefficient becomes positive if and only if
\[
\Delta_{i,t}>\tau_i,\qquad
\tau_i=\log\!\left(1-\frac{A_i}{\beta}\right)>0.
\]

If $b_{i,t}=0$, then
\[
A^{\mathrm{eff}}_{i,t}=A_i<0.
\]
For $b_{i,t}=1$ and $\Delta_{i,t}\leq0$,
\[
A^{\mathrm{eff}}_{i,t}
=A_i+\beta(e^{\Delta_{i,t}}-1)\leq A_i.
\]
For $\Delta_{i,t}>0$, a sign reversal occurs exactly when
\[
\begin{aligned}
A_i+\beta(e^{\Delta_{i,t}}-1)&>0\\
\Longleftrightarrow\quad
e^{\Delta_{i,t}}&>1-\frac{A_i}{\beta}\\
\Longleftrightarrow\quad
\Delta_{i,t}&>
\log\left(1-\frac{A_i}{\beta}\right)=\tau_i.
\end{aligned}
\]
Therefore $A^{\mathrm{eff}}_{i,t}\leq0$ for
$0<\Delta_{i,t}\leq\tau_i$, and it becomes positive only when
$\Delta_{i,t}>\tau_i$.

The same result has a general form. Define
\[
G_{i,t}=\rho_{i,t}A_i<0,
\qquad
\lambda_{i,t}=\beta\widetilde b_{i,t}>0.
\]
Then $G_{i,t}+\lambda_{i,t}(e^{\Delta_{i,t}}-1)$ changes sign at
\begin{equation}
\tau^{\mathrm{gen}}_{i,t}
=\log\!\left(1-\frac{G_{i,t}}{\lambda_{i,t}}\right).
\label{eq:app-general-threshold}
\end{equation}
The threshold in the main paper follows from
$\rho_{i,t}=\widetilde b_{i,t}=1$.

\subsection{Boundary Cases and Scope}
\label{app:theory-boundaries}

If no rollout executes, then $\Rvote=\bot$, $b_{i,t}=0$, and
$\mathcal L_{\mathrm{KL}}=0$. A majority-cluster rollout is also masked out.
If an LP signature matches the reference, or a response stage has no
solver-grounded mapping, the structural gate blocks the contextual update
regardless of the teacher's token probability or any surface-form shift.

The decomposition clarifies the division of labor: the group advantage
provides response-level outcome competition, the teacher ratio provides a
token-level contextual direction, and the LP-derived mask determines where
that direction may act. A negative group signal anchors the update, so a
contextual preference must exceed
Equation~\eqref{eq:app-general-threshold} before reversing its sign.

These results characterize the update mechanism; they are not a performance
guarantee. They do not prove that the majority LP is correct, that every
activated likelihood shift is a semantic correction, or that an optimization
step must improve expected solver accuracy. The derivations require the
fixed-teacher and fixed-mask convention and apply to the unclipped branch;
clamping or clipping can alter the displayed gradient.